\documentclass[12pt,oneside]{article}
\usepackage{amsmath,amssymb,amsfonts,amsthm}
\usepackage{color}
\usepackage[all]{xy}
\usepackage{graphicx}
\usepackage{color}
\usepackage{makeidx}
\usepackage{psfrag}
\usepackage{epstopdf}
\usepackage{multirow}
\usepackage{rotating}
\usepackage{maplestd2e}
\def\paa{\dot{\partial}}

\def\Section#1{\vspace{30truept}\addtocounter{section}{1}\setcounter{thm}{0}\setcounter{equation}{0}
{\noindent\Large\bf\arabic{section}.#1}\par \vspace{12pt}}

\newtheorem{thm}{Theorem}[section]

\newtheorem{rem}[thm]{Remark}

\numberwithin{equation}{section}

\newcommand\undersym[2]{\raisebox{-7pt}{\tiny$#2$}{\kern-8pt}\mbox{$#1$}}
\newcommand\undersymm[2]{\raisebox{-8pt}{\tiny$#2$}{\kern-15pt}\mbox{$#1$}}
\newcommand\overast[1]{\raisebox{10pt}{\small$\ast$}{\kern-7.5pt}\mbox{$#1$}}
\newcommand\overlind[1]{\raisebox{10pt}{\small$\overline{{\hspace{2pt}}\star}$}{\kern-7.5pt}\mbox{$#1$}}
\newcommand\overlinc[1]{\raisebox{10pt}{\tiny$\overline{{\hspace{2pt}}\circ}$}{\kern-7.5pt}\mbox{$#1$}}
\newcommand\overlina[1]{\raisebox{10pt}{\small$\overline{{\hspace{1pt}}\ast}$}{\kern-7.5pt}\mbox{$#1$}}
\newcommand\overcirc[1]{\raisebox{10pt}{\tiny{$\circ$}}{\kern-7.5pt}\mbox{$#1$}}
\newcommand\overdiamond[1]{\raisebox{10pt}{\small$\star$}{\kern-7.5pt}\mbox{$#1$}}

\begin{document}
\title{\bf{New Finsler Package}}
\author{\bf{ Nabil L. Youssef$^{\,1,2}$ and S. G.
Elgendi$^{3}$}}
\date{}
\maketitle                     
\vspace{-1.16cm}
\begin{center}
{$^{1}$Department of Mathematics, Faculty of Science,\\ Cairo
University, Giza, Egypt}\\

\bigskip
$^{2}$Center for Theoretical Physics (CTP)\\
at the Britich University in Egypt (BUE)
\end{center}

\begin{center}
{$^{3}$Department of Mathematics, Faculty of Science,\\ Benha
University, Benha,
 Egypt}
\end{center}

\begin{center}
E-mails: nlyoussef@sci.cu.edu.eg, nlyoussef2003@yahoo.fr\\
{\hspace{1.8cm}}salah.ali@fsci.bu.edu.eg, salahelgendi@yahoo.com
\end{center}
\smallskip
\vspace{1cm} \maketitle
\smallskip
{\vspace{-0.5cm}}
\noindent{\bf \begin{center}
Abstract
\end{center}}
The book \lq \lq Handbook of Finsler geometry"   has been  included with a CD containing  an elegant   Maple package, FINSLER, for  calculations in Finsler geometry. Using this package,     an example concerning  a Finsler generalization of Einstein's vacuum field equations was treated. In this example,  the calculation of the  components of the hv-curvature  of Cartan connection leads to  wrong expressions.  On the other hand,   the FINSLER  package works only  in dimension  four. We introduce a new Finsler package in which we fix the two problems and solve them. Moreover,     we extend this package to compute  not only the geometric objects associated with  Cartan connection but also those associated with  Berwald, Chern and Hashiguchi connections in any  dimension. These  improvements have been illustrated by a concrete  example. Furthermore,  the problem of simplifying  tensor expressions is treated.
This paper is intended  to make calculations in Finsler geometry more easier and simpler.\\


\newpage



\Section{  Introduction}

\par

Antonelli et. al. have had a good contribution in Finsler geometry computations using MAPLE (cf. \cite{Rutz1}, \cite{Rutz2},\cite{Portugal1}).
Rutz and Portugal \cite{hbfinsler} have introduced the remarkable FINSLER package \cite{Rutz3} (it is also included in a CD with the \lq\lq Handbook of Finsler geometry" \cite{hbfinsler1}). They illustrated how to use this package by an example related to general relativity.

\vspace{5pt}
During  the preparation of our paper \cite{ND-cartan}, searching for some Finsler counterexamples, we have encountered some problems concerning FINSLER   package. In fact, we studied  an example in which  the coefficients of Berwald connection are functions of positional argument $x^i$  only. Hence, the  space  under consideration is Berwaldian and is thus Landesbergian. It is well known that for a  Landesberg space the hv-curvature $P^h_{ijk}$ of Cartan connection vanishes. But according to the package, the program calculated non-vanishing components of  $P^h_{ijk}$.
After a deep reading  of the source code (Finsler.mpl), we discovered an  error in the definition of $P^h_{ijk}$  (similar error is found  in  \lq\lq Handbook of Finsler geometry, II \rq\rq, page 1154).
Another problem with this package  is that of dimension. If one considers   a Finsler space of dimension three, the package  can not compute  the components of  the h-curvature $R^h_{ijk}$ and hv-curvature $P^h_{ijk}$ of Cartan connection.

\vspace{5pt}
In our modified package we solve the above two mentioned problems. We illustrate our modification and extension of the FINSLER package by treating a concrete example of a three dimensional Finsler space. We calculate the curvature tensors of the four fundamental connections of Finsler geometry, namely, Cartan, Berwald, Chern and Hashiguchi connections.  The geometric objects, not defined in the FINSLER package, can be added in a similar manner. We also propose a technique for simplifying tensor expressions.


\Section{ Notations and preliminaries}

In this section, we give a brief introduction to Finsler connections.  For more details, we refer, for example, to \cite{hbfinsler1},  \cite{chernshen},  \cite{buctaru} and \cite{ourbook}.

\par Let  $(M,F)$ be a Finsler manifold. Let $(x^i)$ be
the coordinates of any  point of $M$ and $(y^i)$ a
supporting element at this point. Partial differentiation with respect to $x^i$ (resp. $y^i$) will be denoted by $\partial_i$ (resp. $\dot{\partial}_i$).
We use the following notations:
\vspace{5pt}

\noindent $l_i:=\dot{\partial}_iF=g_{ij} l^j=g_{ij}\frac{y^j}{F}$: the
    normalized supporting element; $l^i:=\frac{y^i}{F}$,

\vspace{3pt}\noindent
 $l_{ij}:=\dot{\partial}_il_j$,

\vspace{3pt}\noindent
   $h_{ij}:=Fl_{ij}=g_{ij}-l_il_j$:  the angular metric tensor,

\vspace{3pt}\noindent
   $C_{ijk}:=\frac{1}{2}\dot{\partial}_kg_{ij}=\frac{1}{4}\dot{\partial}_i
    \dot{\partial}_j\dot{\partial}_k F^2$:  the Cartan  tensor,

\vspace{3pt}\noindent
      $C^i_{jk}:=g^{ri}C_{rjk}$:
     the (h)hv-torsion tensor,

\vspace{3pt}\noindent
     $\gamma^i_{jk}(x,y):=\frac{1}{2}g^{ir}(\partial_jg_{kr}+\partial_kg_{jr}-\partial_rg_{jk})$:
        the Christoffel symbols with respect to $\partial_i$,

\vspace{3pt}\noindent
        $G^i(x,y):=\frac{1}{2}\gamma^i_{jk}y^jy^k$: the components of the   canonical spray associated
    with $(M,F)$,

\vspace{3pt}\noindent
    $ N^i_j:=\dot{\partial}_jG^i$: the Barthel (or Cartan nonlinear) connection
    associated with $(M,F)$,

\vspace{3pt}\noindent
    $G^i_{jh}:=\dot{\partial}_hN^i_j=\dot{\partial}_h\dot{\partial}_jG^i$:
     the coefficients of Berwald connection,

\vspace{3pt}\noindent
    $\delta_i:=\partial_i-N^r_i\dot{\partial}_r$: the basis
     vector fields of the horizontal bundle,

\vspace{3pt}\noindent
        $\Gamma^i_{jk}(x,y):=\frac{1}{2}g^{ir}(\delta_jg_{kr}+\delta_kg_{jr}-\delta_rg_{jk})$:
          the Christoffel symbols with respect to $\delta_i$.

\bigskip
A \textbf{Finsler  connection} \cite{Anas} on $M$
is a triple $F\Gamma=( \textbf{F}^{\,i}_{jk}(x,y),
         \textbf{N}^i_j(x,y),\textbf{C}^{\,i}_{jk}(x,y))$ such that,
         under a change of coordinates $(x^i)\rightarrow
         (\widetilde{x}^k)$, the geometric objects $\textbf{F}^{\,i}_{jk}(x,y)$,  $\textbf{N}^i_j(x,y)$ and
         $\textbf{C}^{\,i}_{jk}$ transform respectively as follows:
$$
\widetilde{\textbf{F}}^{\,k}_{ij}=\frac{\partial
\widetilde{x}^k}{\partial x^l}
 \frac{\partial x^p}{\partial \widetilde{x}^i}\frac{\partial x^q}{\partial \widetilde{x}^j}
 \textbf{F}^{\,l}_{pq}+\frac{\partial^2 {x}^p}{\partial \widetilde{x}^i \partial \widetilde{x}^j}
 \frac{\partial \widetilde{x}^k}{\partial {x}^p},$$
$$\widetilde{\textbf{N}}^i_j= \frac{\partial
\widetilde{x}^i}{\partial x^p}\frac{\partial x^q}{\partial
\widetilde{x}^j}\textbf{N}^p_q+\frac{\partial {x}^p}{\partial
\widetilde{x}^j}\frac{\partial^2 \widetilde{x}^i}{\partial x^p
\partial x^q}y^q,\quad
\widetilde{\textbf{C}}^{\,k}_{ij}=\frac{\partial
\widetilde{x}^k}{\partial x^l}
 \frac{\partial x^p}{\partial \widetilde{x}^i}\frac{\partial x^q}{\partial
 \widetilde{x}^j}\textbf{C}^{\,l}_{pq}.$$
 Moreover, $F\Gamma$ defines  two types of
covariant derivatives:
\begin{eqnarray*}
  X^i_{j\mid k} &:=& \delta_kX^i_j+X^r_j\textbf{F}^{\,i}_{rk}-X^i_r\textbf{F}^{\,r}_{jk}. \\
  X^i_j|_k      &:=& \dot{\partial}_kX^i_j+X^r_j\textbf{C}^{\,i}_{rk}-X^i_r\textbf{C}^{\,r}_{jk}.
\end{eqnarray*}

Let $F\Gamma=( \textbf{F}^{\,i}_{jk},
         \textbf{N}^i_j,\textbf{C}^{\,i}_{jk})$ be an arbitrary Finsler connection.
The  (h)h-,  (h)hv-,  (v)h-,  (v)hv- and (v)v-torsion tensors of
$F\Gamma$ are  given respectively by~\cite{r2.6}:
$$
 \textbf{T}^i_{jk}=\textbf{F}^i_{jk}-\textbf{F}^i_{kj},\qquad
 \textbf{C}^i_{jk}= \mathrm{the \ connection \  parameters} \  \textbf{C}^i_{jk},$$
 $$
   \textbf{R}^i_{jk}=\delta_k\textbf{N}^i_j-\delta_j\textbf{N}^i_k,\qquad
    \textbf{P}^i_{jk}=\dot{\partial}_k\textbf{N}^i_j-\textbf{F}^i_{jk},\qquad
    \textbf{S}^i_{jk}=\textbf{C}^i_{jk}-\textbf{C}^i_{kj}.
$$
and the h-,  hv- and v-curvature  tensors of $F\Gamma$ are
 given respectively by \cite{r2.6}:
$$
  \textbf{R}^i_{hjk} = \mathfrak{A}_{(j,k)}\{{\delta_k\textbf{F}^i_{hj}}
+\textbf{F}^m_{hj}\textbf{F}^i_{mk}\}+\textbf{C}^i_{hm}\textbf{R}^m_{jk}, $$
  $$\textbf{ P}^i_{hjk}=\dot{\partial}_k\textbf{F}^i_{hj}-\textbf{C}^i_{hk\mid
j}+\textbf{C}^i_{hm}\textbf{P}^m_{jk},\qquad
  \textbf{S}^i_{hjk}=\mathfrak{A}_{(j,k)}\{\dot{\partial}_k\textbf{C}^i_{hj}+\textbf{C}^m_{hk}\textbf{C}^i_{mj}\},
$$
where $\mathfrak{A}_{(j,k)}\{A_{jk}\}:=A_{jk}-A_{kj}$.\\

The \textbf{Cartan connection} is given by
$C\Gamma=(\Gamma^i_{jk},N^i_j,C^i_{jk})$,
where $\Gamma^i_{jk}$, $N^i_j$ and ${C}^i_{jk}$ are as defined  above.
The (h)hv-,  (v)h- and (v)hv-torsion tensors of $C\Gamma$ are:\vspace{-5pt}
 $$C^i_{jk}=\frac{1}{2}g^{ir}\dot{\partial}_kg_{rj}, \qquad R^i_{jk}=\delta_kN^i_j-\delta_jN^i_k,\qquad P^i_{jk}=\dot{\partial}_kN^i_j-\Gamma^i_{jk}. $$
The h-,  hv- and v-curvature  tensors of $C\Gamma$ are:

  $$R^i_{hjk} =\mathfrak{A}_{(j,k)}\{{\delta_k\Gamma^i_{hj}}
+\Gamma^m_{hj}\Gamma^i_{mk}\}+C^i_{hm}R^m_{jk}, $$
   $$P^i_{hjk}= \dot{\partial}_k\Gamma^i_{hj}-C^i_{hk\mid
j}+C^i_{hm}P^m_{jk},\quad
  S^i_{hjk}=\mathfrak{A}_{(j,k)}\{C^m_{hk}C^i_{mj}\}.$$


The \textbf{Berwald  connection}  is given by
$B\Gamma=(G^i_{jk},N^i_j,0)$. The associated geometric objects will be
marked by a circle. The (v)h-torsion tensor of $B\Gamma$ is  given by:
$$\overcirc{R}^i_{jk}={R}^i_{jk}=\delta_kN^i_j-\delta_jN^i_k.$$
\noindent The h-,  and hv-curvature  tensors of $B\Gamma$ are:
 $$ \overcirc{R}^i_{hjk} = \mathfrak{A}_{(j,k)}\{{\delta_kG^i_{hj}}
+G^m_{hj}G^i_{mk}\},\qquad
   \overcirc{P}^i_{hjk}=\dot{\partial}_kG^i_{hj}.$$


 The  \textbf{Chern (Rund) connection} is  given by $R\Gamma=(\Gamma^i_{jk},N^i_j,0)$.
The associated geometric objects will be
marked by a star. The (v)h- and (v)hv-torsion tensors of $R\Gamma$ are:
  $$\overdiamond{R}^i_{jk}={R}^i_{jk}=\delta_kN^i_j-\delta_jN^i_k,\qquad
  \overdiamond{P}^i_{jk}=P^i_{jk}=\dot{\partial}_kN^i_j-\Gamma^i_{jk}.$$
\noindent The h- and  hv-curvature  tensors of $R\Gamma$ are:
 $$ \overdiamond{R}^i_{hjk} = \mathfrak{A}_{(j,k)}\{{\delta_k\Gamma^i_{hj}}
+\Gamma^m_{hj}\Gamma^i_{mk}\}, \qquad
   \overdiamond{P}^i_{hjk}=\dot{\partial}_k\Gamma^i_{hj}.$$


 The  \textbf{Hashiguchi  connection} is given by
$H\Gamma=(G^i_{jk},N^i_j,C^i_{jk})$.
The associated geometric objects will be marked by an asterisk. The (h)hv-  and  (v)h-torsion tensors of~$H\Gamma$~are:
 $$\overast{C}^i_{jk}=  C^i_{jk} ,\qquad
   \overast{R}^i_{jk}={R}^i_{jk}=\delta_kN^i_j-\delta_jN^i_k.$$
\noindent The h-,  hv- and v-curvature  tensors of $H\Gamma$ are:
  $$\overast{R}^i_{hjk} =\mathfrak{A}_{(j,k)}\{{\delta_kG^i_{hj}}
+G^m_{hj}G^i_{mk}\}+C^i_{hm}R^m_{jk}, $$
   $$\overast{P}^i_{hjk}= \dot{\partial}_kG^i_{hj}-C^i_{hk\stackrel{*}|
j},\qquad
  \overast{S}^i_{hjk}=\mathfrak{A}_{(j,k)}\{C^m_{hk}C^i_{mj}\}.$$

\begin{center}{{Table 1: Fundamental linear connections \cite{ourbook}}}
\\[.4 cm]

\small{\begin{tabular} {|c|c|c|c|c|c|}\hline
 &\multirow{2}{*}{\textbf{} }&\multirow{2}{*}{{ Cartan }}&\multirow{2}{*}{{
 Berwald}}&
 \multirow{2}{*}{ { Chern (Rund)  }} &\multirow{2}{*}{{ Hashiguchi }}\\
 &&&&&\\[0.1 cm]\cline{2-6}
&\multirow{2}{*}{$(\textbf{F}^{h}_{ij},
\textbf{N}^{h}_{i},\textbf{C}^{h}_{ij})$}&
\multirow{2}{*}{$(\Gamma^{h}_{ij}, N^{h}_{i}, C^{h}_{ij})$}&\multirow{2}{*}{$( G^{h}_{ij}, N^{h}_{i}, 0)$}&
\multirow{2}{*}{$( \Gamma^{h}_{ij}, N^{h}_{i}, 0)$}&
\multirow{2}{*}{$( G^{h}_{ij}, N^{h}_{i}, C^{h}_{ij})$}\\
\rule{.2cm}{0pt}
\begin{rotate}{90}\hspace{-.15cm}Connection\end{rotate}\rule{.3cm}{0pt}&&
& & &
\\[0.1 cm]\hline
&\multirow{2}{*}{{ (h)h-torsion} $\textbf{T}^{i}_{jk}$}&\multirow{2}{*}{$0$}&\multirow{2}{*}{$0$}&\multirow{2}{*}{$0$}&\multirow{2}{*}{$0$}\\
&\multirow{2}{*}{{ (h)hv-torsion}
$\textbf{C}^{i}_{jk}$}&\multirow{2}{*}{$C^{i}_{jk}$} &\multirow{2}{*}{$0$}&
\multirow{2}{*}{$0$}&
\multirow{2}{*}{$C^{i}_{jk}$}
\\[0.1 cm]&&&&&
\\[0.1 cm]\cline{2-6}
&&&&&\\
\rule{.3cm}{0pt}
\begin{rotate}{90}Torsions\end{rotate}\rule{.3cm}{0pt}
&{ (v)h-torsion} $\textbf{R}^{i}_{jk}$&${R}^{i}_{jk}$&$\overcirc{R}^{i}_{jk}={R}^{i}_{jk}$
&$\overdiamond{R}^{i}_{jk}={R}^{i}_{jk}$&
$\overast{R}^{i}_{jk}=R^{i}_{jk}$
\\[0.1 cm]
&{ (v)hv-torsion}
$\textbf{P}^{i}_{jk}$&${P}^{i}_{jk}=C^{i}_{jk|h}y^h$&$0$&
  $\overdiamond{P}^{i}_{jk}=P^{i}_{jk}$&$0$
\\[0.1 cm]
&{ (v)v-torsion} $\textbf{S}^{i}_{jk}$&$0$ &$0$&$0$&$0$
\\[0.1 cm]\hline
&&&&&\\
&{ h-curvature} $\textbf{R}^{h}_{ijk}$& ${R}^{h}_{ijk}$&$\overcirc{R}^{h}_{ijk}$&
$\overdiamond{R}^{h}_{ijk}$&
$\overast{R}^{h}_{ijk}$
\\[0.1 cm]
&{ hv-curvature} $\textbf{P}^{h}_{ijk}$& ${P}^{h}_{ijk}$&$\overcirc{P}^{h}_{ijk}$&
$\overdiamond{P}^{h}_{ijk}$&
$\overast{P}^{h}_{ijk}$
\\[0.1 cm]
\rule{.2cm}{0pt}
\begin{rotate}{90}\hspace{.2cm}Curvatures\end{rotate}\rule{.2cm}{0pt}&{ v-curvature}
$\textbf{S}^{h}_{ijk}$& ${S}^{h}_{ijk}$&$0$&
$0$&$\overast{S}^{h}_{ijk}=S^{h}_{ijk}$
\\[0.1 cm]\hline &&&&&\\
&{ h-cov. der.}&
${K}^{i}_{j|k}$&${K}^{i}_{j\stackrel{\circ}|k}$&${K}^{i}_{j\stackrel{\star}|k}={K}^{i}_{j|k}$ &
${K}^{i}_{j\stackrel{*}|k}={K}^{i}_{j\stackrel{\circ}|k}$
\\{\vspace{-7pt}}&&&&&\\
&{ v-cov. der.}&
${K}^{i}_{j}|_k$&${K}^{i}_{j}{\stackrel{\circ}|}_k=\paa_{k}{K}^{i}_{j}$&${K}^{i}_{j}{\stackrel{\star}|}_k={K}^{i}_{j}{\stackrel{\circ}|}_k$
&${K}^{i}_{j}{\stackrel{*}|}_k={K}^{i}_{j}|_k$

\\
\rule{.3cm}{0pt}
\begin{rotate}{90}\hspace{.2cm} Covariant \end{rotate}\rule{.3cm}{0pt}
\rule{0cm}{0pt}\begin{rotate}{90}
\hspace{.3cm}derivatives\end{rotate}\rule{.3cm}{0pt}&&&&&\\ \hline
\end{tabular}}
\end{center}

{\vspace{4pt}}
\Section{ Notes on  the FINSLER package }

    In \cite{hbfinsler}, Rutz and Portugal  discussed and  applied  the FINSLER package they introduced in \cite{Rutz3}. This package is an  extension of  the RIEMANN package \cite{Portugal1}. The FINSLER package is included in a CD with the book \cite{hbfinsler1}, where an interesting example related to general relativity, namely, a family of metrics known as the Schwarzschild solution to Einstein's field equations, has been treated. Important geometric objects and, in particular, the three curvature tensors of Cartan connection have been computed.

\smallskip
When  performing some applications using the FINSLER package, we have encountered  some problems. To show one of these problems, let us  consider  the following example.
\smallskip\noindent
Let $M =\mathbb{R}^4$,
$U=\{(x,y)\in\mathbb{R}^4\times \mathbb{R}^4:x_1\neq 0;\,y_4\neq 0,\, y_1^2+y_2^2+y_3^2\neq 0\}$.
Let $F$ be the Finsler structure  defined on the open subset $U$ of $TM$ by:
  $$F=\sqrt{x_1 y_4 \sqrt{y_1^2+y_2^2+y_3^2}}.$$
 Based on this package, the non-vanishing coefficients of  Berwald connection are  as follows:
  $$G^{1}_{11}=G^{2}_{12}=G^{3}_{13}=\frac{1}{x_1}, \quad G^{1}_{22}=G^{1}_{33}=-\frac{1}{x_1}.$$
This shows that the coefficients of Berwald connection are functions of the positional  argument $x^i$  only. Hence, the  space  under consideration is Berwaldian and is thus Landesbergian. Consequently,  the hv-curvature $P^h_{ijk}$ of Cartan connection should vanish identically. However,  the FINSLER package calculated non-vanishing components of  $P^h_{ijk}$.

\smallskip
After a deep study  of the source code (Finsler.mpl), we have discovered some wrong indices  in the definition of $P^h_{ijk}$. (Similar error is found  in  \cite{hbfinsler1},  page 1154).
Another problem with this package  is the problem of dimension. If one considers   a three dimensional Finsler space, the package  can not compute  the components of  the hh-curvature $R^h_{ijk}$ and hv-curvature $P^h_{ijk}$ of Cartan connection. The package response is that these objects are outside dimension.

\smallskip
Summing up, we have two problems with the Rutz and Portugal's package. The first  is the wrong calculations of the  curvature $P^h_{ijk}$. The second  is the disability of computing  $R^h_{ijk}$ and $P^h_{ijk}$  in dimensions different from $4$.


\Section{ Improvement of the package }

In this section, we  solve  the two above   mentioned  problems. Moreover, we  extend   the package in order   to compute various  geometric objects associated not only with Cartan connection but also with  the other fundamental  connections in Finsler geometry.  And this is for any dimension. Other   geometric objects can be similarly added    to the package. We illustrate these tasks using  a concrete example.

\smallskip
 Rutz and Portugal have  illustrated how to use the  package \cite{hbfinsler}. However, let  us  recall  some
 instructions  to make  the use of this package easier.  When we write, for example, \linebreak N[i,-j] we  mean $N^i_j$, i.e., a positive (resp. negative) index means that it is a contravariant (resp. covariant) index. If one wants to  lower or raise  an index by  the metric or the inverse  metric, he just changes its sign from  positive to negative or  vice versa.   The command  \emph{tdiff}(N[i,-j], X[k]) means ${\partial}_k N^i_j$, the command  \emph{tddiff}(N[i,-j], Y[k]) means $\dot{\partial}_k N^i_j$ and the command \emph{Hdiff}(N[i,-j], X[k]) means $\delta_k N^i_j$.

\smallskip
In addition to the definitions  of geometric objects  existing  already in the FINSLER package, we add other definitions by using the command \emph{definetensor}.  We rewrite the correct expression of $P^h_{ijk}$  and tackle the issue of dimension.

\smallskip
Now, let us illustrate what have   been said before  using a concrete example.

\smallskip
\noindent Let $M=\mathbb{R}^3$, $U=\{(x1,x2,x3;y1,y2,y3)\in \mathbb{R}^3 \times \mathbb{R}^3: x3\neq 0;\, y2\neq 0, \,  y1^2+y3^2\neq 0 \}\subset TM$.  Let  $F$ be the Finsler structure  defined on $U$ by
$$F=\sqrt{\frac{x3y1^3}{y2}+y3^2}.$$
It should first be noted that, according to Table 1,  we have only three  independent torsions, namely, $C^h_{ij}$, $R^h_{ij}$ and $P^h_{ij}$. So, we will compute these torsions for  Cartan connection  and we will not repeat their calculation for the  other connections.

\smallskip
Following the instructions of  the  FINSLER package,    the following  calculations can be  performed.
\bigskip

\begin{maplegroup}
\begin{mapleinput}
\mapleinline{active}{1d}{restart;
libname := libname, `c:/Finsler`:
with(Finsler);}{}
\end{mapleinput}
\mapleresult
\begin{maplelatex}
\mapleinline{inert}{2d}{[Dcoordinates, Hdiff, K, connection, init, metricfunction, tddiff]}{\[ [{\it Dcoordinates},{\it Hdiff},K,{\it connection},{\it init},{\it metricfunction},{\it tddiff}]\]}
\end{maplelatex}
\end{maplegroup}
\begin{maplegroup}
\begin{mapleinput}
\mapleinline{active}{1d}{Dimension := 3:
coordinates(x1,x2,x3):
Dcoordinates(y1,y2,y3):}{}
\end{mapleinput}
\mapleresult
\begin{maplelatex}
\mapleinline{inert}{2d}{`The coordinates are:`}{\[ \mbox {{\tt `The coordinates are:`}}\]}
\end{maplelatex}
\mapleresult
\begin{maplelatex}
\mapleinline{inert}{2d}{X*` `^`1` = x1}{\[ X^1={\it x1}\]}
\end{maplelatex}
\mapleresult
\begin{maplelatex}
\mapleinline{inert}{2d}{X*` `^2 = x2}{\[ X^{2}={\it x2}\]}
\end{maplelatex}
\mapleresult
\begin{maplelatex}
\mapleinline{inert}{2d}{X*` `^3 = x3}{\[ X^{3}={\it x3}\]}
\end{maplelatex}
\mapleresult
\mapleresult
\begin{maplelatex}
\mapleinline{inert}{2d}{`The d-coordinates are:`}{\[ \mbox {{\tt `The d-coordinates are:`}}\]}
\end{maplelatex}
\mapleresult
\begin{maplelatex}
\mapleinline{inert}{2d}{Y*` `^`1` = y1}{\[ Y^{\mbox {{1}}}={\it y1}\]}
\end{maplelatex}
\mapleresult
\begin{maplelatex}
\mapleinline{inert}{2d}{Y*` `^2 = y2}{\[ Y^{2}={\it y2}\]}
\end{maplelatex}
\mapleresult
\begin{maplelatex}
\mapleinline{inert}{2d}{Y*` `^3 = y3}{\[ Y^{3}={\it y3}\]}
\end{maplelatex}
\end{maplegroup}
\vspace{.2cm}
\begin{maplegroup}
\begin{Maple Normal}{
\textbf{Finsler structure F:}}\end{Maple Normal}

\end{maplegroup}
\vspace{.2cm}

\begin{maplegroup}
\begin{mapleinput}
\mapleinline{active}{1d}{F := sqrt(x3*y1\symbol{94}3/y2+y3\symbol{94}2); }{\[\]}
\end{mapleinput}
\mapleresult
\begin{maplelatex}
\mapleinline{inert}{2d}{F := sqrt(y1^3*x3/y2+y3^2)}{\[ F\, := \, \sqrt{{\frac {{{\it x3}{\it y1}}^{3}}{{\it y2}}}+{{\it y3}}^{2}}\]}
\end{maplelatex}
\end{maplegroup}
\vspace{.2cm}

\begin{maplegroup}
\begin{Maple Normal}{
\textbf{Plotting  the Finsler structure in a special domain}\textbf{:}}\end{Maple Normal}

\end{maplegroup}
\vspace{.2cm}
\begin{maplegroup}
\begin{mapleinput}
\mapleinline{active}{1d}{plot3d(subs(x1=5,y3=5,F), y1 = -2..2, y2 =-2..2,
axes=BOXED,style=patch); }{\[\]}
\end{mapleinput}
\mapleresult
\begin{center}
\includegraphics[width=1\textwidth]{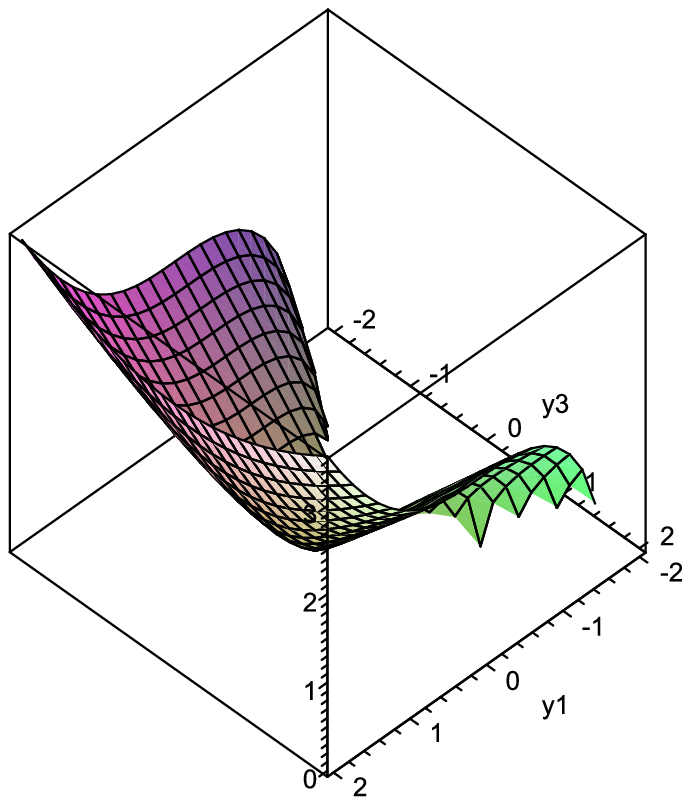}
\end{center}
\vspace{-1.1cm}
\begin{center}
     Figure 1
\end{center}

\end{maplegroup}

\vspace{.2cm}

\begin{maplegroup}
\begin{Maple Normal}{
\textbf{Metric tensor $g_{ij}$:} }\end{Maple Normal}

\end{maplegroup}

\vspace{.2cm}
\begin{maplegroup}
\begin{mapleinput}
\mapleinline{active}{1d}{F0:= y1^3*x3/y2+y3^2;
}{}
\end{mapleinput}
\mapleresult
\begin{maplelatex}
\mapleinline{inert}{2d}{F0 := y1^3*x3/y2+y3^2}{\[ {\it F0}\, := \,{\frac {{\it x3}{{\it y1}}^{3}}{{\it y2}}}+{{\it y3}}^{2}\]}
\end{maplelatex}
\end{maplegroup}
\begin{maplegroup}
\begin{mapleinput}
\mapleinline{active}{1d}{metricfunction(F0):

}{}
\end{mapleinput}
\mapleresult
\begin{maplelatex}
\mapleinline{inert}{2d}{`The components of the metric are:`}{\[ \mbox {{\tt `The components of the metric are:`}}\]}
\end{maplelatex}
\mapleresult
\end{maplegroup}
{\vspace{-.5cm}}\begin{eqnarray*}
  g_{{\it x1 x1} }&=&\frac{3~{\it x3}{\it y1} }{{\it y2} }\quad\quad\quad  g_{{\it x1 x2} }=-\frac{3}{2}~\frac{{\it x3}{\it y1}^{2} }{{\it y2}^{2}}\\
  g_{{\it x2 x2} }&=&\frac{{\it x3}{\it y1}^{3} }{{\it y2}^{3}}\quad\quad\quad\,\,  g_{{\it x3 x3} }=1
\end{eqnarray*}
\begin{maplegroup}
\begin{Maple Normal}{
\textbf{Inverse  metric tensor $g^{ij}$:} }\end{Maple Normal}

\end{maplegroup}
\vspace{.2cm}
\begin{maplegroup}
\begin{mapleinput}
\mapleinline{active}{1d}{show(g[i,j]); }{\[\]}
\end{mapleinput}
\end{maplegroup}
{\vspace{-.7cm}}\begin{eqnarray*}
g^{{\it x1 x1}}&=&{\frac {{4\it  y2}}{{3\it x3}{\it y1}}}\quad\quad\quad  g^{{\it x1 x2}}={\frac {{{2\it y2}}^{2}}{{{\it x3}{\it y1}}^{2}}}\\
g^{{\it x2 x2}}&=&{\frac {{{4\it  y2}}^{3}}{{{\it x3}{\it y1}}^{3}}}\quad\quad\quad  g^{{\it x3 x3}}=1
\end{eqnarray*}
\begin{maplegroup}
\begin{Maple Normal}{
\textbf{Supporting element $l_i$:}}\end{Maple Normal}

\end{maplegroup}

\vspace{.2cm}

\begin{maplegroup}
\begin{mapleinput}
\mapleinline{active}{1d}{show(l[-i]);
}{}
\end{mapleinput}
\mapleresult
\end{maplegroup}
\begin{eqnarray*}
 l _{{\it x1} }&=&\frac{3}{2}\frac{{\it x3}{\it y1}^{2} }{{\it y2} \sqrt{\frac{{\it x3}{\it y1}^{3} +{\it y2}{\it y3}^{2} }{{\it y2} }}}\quad\quad\quad  l _{{\it x2} }=-\frac{1}{2}\frac{{\it x3}{\it y1}^{3} }{{\it y2}^{2}\sqrt{\frac{{\it x3}{\it y1}^{3} +{\it y2}{\it y3}^{2} }{{\it y2} }}}\\
 l _{{\it x3} }&=&\frac{{\it y3} }{\sqrt{\frac{{\it x3}{\it y1}^{3} +{\it y2}{\it y3}^{2} }{{\it y2} }}}
\end{eqnarray*}
\begin{maplegroup}
\begin{Maple Normal}{
\textbf{Angular metric tensor $h_{ij}$:}}\end{Maple Normal}

\end{maplegroup}

\vspace{.2cm}

\begin{maplegroup}
\begin{mapleinput}
\mapleinline{active}{1d}{definetensor(h[-i,-j] = g[-i,-j]-l[-i]*l[-j],symm);
}{}
\end{mapleinput}
\mapleresult
\end{maplegroup}
\begin{maplegroup}
\begin{mapleinput}
\mapleinline{active}{1d}{show(h[-i, -j]); }{\[\]}
\end{mapleinput}
\end{maplegroup}
{\vspace{-.5cm}}\begin{eqnarray*}
h _{x1x1}&=&\frac{3}{4}\frac{{\it x3}{\it y1} \left({\it x3}{\it y1}^{3} +4{\it y2}{\it y3}^{2} \right)}{{\it y2} \left({\it x3}{\it y1}^{3} +{\it y2}{\it y3}^{2} \right)}\quad\quad\quad  h _{x1x2}=-\frac{3}{4}\frac{{\it x3}{\it y1}^{2} \left({\it x3}{\it y1}^{3} +2{\it y2}{\it y3}^{2} \right)}{{\it y2}\left({\it x3}{\it y1}^{3} +{\it y2}{\it y3}^{2} \right)}\\
h _{x_1x_3 }&=&-\frac{3}{2}\frac{{\it x3}{\it y1}^{2}{\it y3} }{{\it x3}{\it y1}^{3} +{\it y2}{\it y3}^{2} }\hspace{2.5cm} h _{x_2x_2 }=\frac{1}{4}\frac{{\it x3}{\it y1}^{3} \left(3{\it x3}{\it y1}^{3} +4{\it y2}{\it y3}^{2} \right)}{{\it y2}^{3}\left({\it x3}{\it y1}^{3} +{\it y2}{\it y3}^{2} \right)}\\
\hspace{-2cm} h _{x_2x_3 }&=&\frac{1}{2}\frac{{\it x3}{\it y1}^{3}{\it y3} }{{\it y2} \left({\it x3}{\it y1}^{3} +{\it y2}{\it y3}^{2} \right)}\quad\quad\quad \,\,\quad h _{x_3x_3 }=\frac{{\it x3}{\it y1}^{3} }{{\it x3}{\it y1}^{3} +{\it y2}{\it y3}^{2} }
\end{eqnarray*}
\begin{maplegroup}
\begin{Maple Normal}{
\textbf{Cartan tensor $C_{ijk}$:}}\end{Maple Normal}

\end{maplegroup}

\vspace{.2cm}

\begin{maplegroup}
\begin{mapleinput}
\mapleinline{active}{1d}{show(C[-i,-j,-k]); }{\[\]}
\end{mapleinput}
\mapleresult
\end{maplegroup}
{\vspace{-.5cm}}\begin{eqnarray*}
  C _{x1x1x_1} &=&\frac{3}{2}\frac{{\it x3} }{{\it y2} }\quad\quad\quad \,\, C _{x_1x1x2 }=-\frac{3}{2}\frac{{\it x3}{\it y1} }{{\it y2}^{2}} \\
  C _{x1x2x_2}&=&\frac{3}{2}\frac{{\it x3}{\it y1}^{2} }{{\it y2}^{3}}\quad\quad  C _{x_2x_2x_2 }=-\frac{3}{2}\frac{{\it x3}{\it y1}^{3} }{{\it y2}^{4}}
\end{eqnarray*}
\begin{maplegroup}
\begin{Maple Normal}{
\textbf{Spray coefficients $G^{i}$:}}\end{Maple Normal}

\end{maplegroup}
\begin{maplegroup}
\begin{mapleinput}
\mapleinline{active}{1d}{show(G[i]);
}{}
\end{mapleinput}
\end{maplegroup}
{\vspace{-.5cm}}\begin{eqnarray*}
 G^{{\it x1}}&=&\frac{1}{2}{\frac {{\it y1}\,{\it y3}}{{\it x3}}}\quad\quad\quad  G^{{\it x2}}=\frac{1}{2}\,{\frac {{\it y2}\,{\it y3}}{{\it x3}}}\\
 G^{{\it x3}}&=&-\frac{1}{4}\,{\frac {{{\it y1}}^{3}}{{\it y2}}}
\end{eqnarray*}
\begin{maplegroup}
\begin{Maple Normal}{
\textbf{Nonlinear connection (Barthel connection) $N^{i}_{j}$:}}\end{Maple Normal}
\end{maplegroup}
\begin{maplegroup}
\begin{mapleinput}
\mapleinline{active}{1d}{show(N[i,-j]);
}{}
\end{mapleinput}
\mapleresult
\end{maplegroup}
{\vspace{-.5cm}}\begin{eqnarray*}
 N^{x1}_{x1} &=&\frac{1}{2}\frac{{\it y3} }{{\it x3} }\quad\quad\quad  N^{x1}_{x3}=\frac{1}{2}\frac{{\it y1} }{{\it x3} }\\
 N^{x2}_{x2}&=&\frac{1}{2}\frac{{\it y3} }{{\it x3} }\quad\quad\quad  N^{x2}_{x3}=\frac{1}{2}\frac{{\it y2} }{{\it x3} }\\
 N^{x3}_{x1}&=&-\frac{3}{4}\frac{{\it y1}^{2}}{{\it y2} }\quad\,\,\,\,  N^{x3}_{x2}=\frac{1}{4}\frac{{\it y1}^{3}}{{\it y2}}
\end{eqnarray*}
\begin{maplegroup}
\begin{Maple Normal}{
\textbf{Coefficients of Berwald connection $G^i_{jk}$:}}\end{Maple Normal}
\end{maplegroup}

\vspace{.2cm}

\begin{maplegroup}
\begin{mapleinput}
\mapleinline{active}{1d}{show(G[i,-j,-k]);
}{}
\end{mapleinput}
\end{maplegroup}
\begin{Maple Normal}{
\begin{Maple Normal}{
}\end{Maple Normal}
}\end{Maple Normal}
{\vspace{-.5cm}}\begin{eqnarray*}
  G^{{\it x3} }_{x1x1}&=&-\frac{3}{2}\frac{{\it y1} }{{\it y2} }\quad\quad\quad  G^{{\it x3} }_{x1x2}=\frac{3}{4}\frac{{\it y1}^{2}}{{\it y2}^{2}} \\
  G^{{\it x1} }_{x1x3}&=&\frac{1}{2{\it x3} }\quad\quad\quad\quad  G^{x3}_{x2x2}=-\frac{1}{2}\frac{{\it y1}^{3}}{{\it y2}^{3}}\\
  G^{x2}_{x2x3}&=&\frac{1}{2{\it x3} }
\end{eqnarray*}
\begin{maplegroup}
\begin{Maple Normal}{
\textbf{Coefficients of Cartan connection $\Gamma^i_{jk}$:}}\end{Maple Normal}
\end{maplegroup}

\vspace{.2cm}

\begin{maplegroup}
\begin{mapleinput}
\mapleinline{active}{1d}{show(Gammastar[i,-j,-k]); }{\[\]}
\end{mapleinput}
\mapleresult
\end{maplegroup}
{\vspace{-.5cm}}\begin{eqnarray*}
  {\it Gammastar}^{{\it x1} }_{x1x1}&=&\frac{1}{2}\frac{{\it y3} }{{\it x3}{\it y1} }\quad\quad\quad  {\it Gammastar}^{{\it x2} }_{x1x1}=\frac{3}{2}\frac{{\it y2}{\it y3} }{{\it x3}{\it y1}^{2} }  \\
   {\it Gammastar}^{{\it x3} }_{x1x1}&=&-\frac{3}{2}\frac{{\it y1} }{{\it y2} }\quad\quad\quad\,\,  {\it Gammastar}^{{\it x1} }_{x1x2}=-\frac{1}{2}\frac{{\it y3} }{{\it x3}{\it y2} }\\
   {\it Gammastar}^{x2}_{x1x2}&=&-\frac{3}{2}\frac{{\it y3} }{{\it x3}{\it y1} }\quad\quad\,\,  {\it Gammastar}^{x3}_{x1x2 }=\frac{3}{4}\frac{{\it y1}^{2}}{{\it y2}^{2}}\\
   {\it Gammastar}^{{\it x1} }_{x1x3}&=&\frac{1}{2{\it x3} }\hspace{1.8cm}  {\it Gammastar}^{{\it x1} }_{x2x2}=\frac{1}{2}\frac{{\it y1}{\it y3} }{{\it x3}{\it y2}^{2}}\\
   {\it Gammastar}^{{\it x2} }_{x2x2}&=&\frac{3}{2}\frac{{\it y3} }{{\it x3}{\it y2} }\hspace{1.3cm}  {\it Gammastar}^{x3}_{x2x2}=-\frac{1}{2}\frac{{\it y1}^{3}}{{\it y2}^{3}}\\
   {\it Gammastar}^{x2}_{x2x3 }&=&\frac{1}{2{\it x3} }
\end{eqnarray*}
\begin{maplegroup}
\begin{Maple Normal}{
\textbf{Torsion tensors of Cartan connection}}\end{Maple Normal}
\end{maplegroup}

\vspace{.2cm}

\begin{maplegroup}
\begin{Maple Normal}{
$\bullet$ \textbf{(h)hv-torsion $C^h_{ij}$:}}\end{Maple Normal}
\end{maplegroup}
\begin{maplegroup}
\begin{mapleinput}
\mapleinline{active}{1d}{show(C[i,-j,-k]);}{\[\]}
\end{mapleinput}
\mapleresult
\end{maplegroup}
{\vspace{-.5cm}}\begin{eqnarray*}
   {\it C}^{{\it x1} }_{x1x1}&=& -\frac{1}{{\it y1}} {\hspace{1.9cm}} {\it C}^{x1}_{x1x2} =\frac{1}{{\it y2} }\\
   {\it C}^{{\it x1} }_{x2x2}&=&-\frac{{\it y1} }{{\it y2}^{2}} \quad\quad\quad   {\it C}^x2_{x1x1} =-\frac{3{\it y2} }{{\it y1}^{2}}\\
  {\it C}^{{\it x2} }_{x1x2}&=&\frac{3}{{\it y1} } {\hspace{2.2cm}}{\it C}^{x2}_{x2x2} =-\frac{3}{{\it y2} }
\end{eqnarray*}
\begin{maplegroup}
\begin{Maple Normal}{
$\bullet$ \textbf{(v)h-torsion $R^h_{ij}$:}}\end{Maple Normal}
\end{maplegroup}
\begin{maplegroup}
\begin{mapleinput}
\mapleinline{active}{1d}{definetensor(RN[i,-j,-k]=Hdiff(N[i,-j],X[k])-Hdiff(N[i,-k],X[j]));}{\[\]}
\end{mapleinput}
\end{maplegroup}
\begin{maplegroup}
\begin{mapleinput}
\mapleinline{active}{1d}{show(RN[i,-j,-k]);}{\[\]}
\end{mapleinput}
\mapleresult
\end{maplegroup}
{\vspace{-.5cm}}\begin{eqnarray*}
  {\it RN}^{{\it x1} }_{x1x2}&=&-\frac{1}{8}\frac{{\it y1}^{3}}{{\it x3}{\it y2}^{2} }\quad\quad  {\it RN}^{x2}_{x1x2} =-\frac{3}{8}\frac{{\it y1}^{2}}{{\it x3}{\it y2} } \\
  {\it RN}^{{\it x1} }_{x1x3}&=&-\frac{1}{4}\frac{{\it y3} }{{\it x3}^{2}}\quad\quad\quad  {\it RN}^{x3}_{x1x3} =\frac{3}{8}\frac{{\it y1}^{2}}{{\it x3}{\it y2} }\\
  {\it RN}^{{\it x2} }_{x2x3}&=&-\frac{1}{4}\frac{{\it y3} }{{\it x3}^{2}}\quad\quad\quad  {\it RN}^{x3}_{x2x3} =-\frac{1}{8}\frac{{\it y1}^{3}}{{\it x3}{\it y2}^{2} }
\end{eqnarray*}
\begin{maplegroup}
\begin{Maple Normal}{
$\bullet$ \textbf{ (v)hv-torsion $P^h_{ij}$:}}\end{Maple Normal}
\end{maplegroup}
\begin{maplegroup}
\begin{mapleinput}
\mapleinline{active}{1d}{definetensor(PT[i,-j,-k] = G[i,-j,-k]- Gammastar[i,-j,-k]):}{\[\]}
\end{mapleinput}
\end{maplegroup}
\begin{maplegroup}
\begin{mapleinput}
\mapleinline{active}{1d}{show(PT[i,-j,-k]); }{\[\]}
\end{mapleinput}
\end{maplegroup}
\mapleresult
{\vspace{-.5cm}}\begin{eqnarray*}
  {\it PT^{x1}_{x1x1} }&=&-\frac{1}{2}\frac{{\it y3} }{{\it x3}{\it y1} }\quad\quad\quad {\it PT^{x2}_{x1x1} }=-\frac{3}{2}\frac{{\it y2}{\it y3} }{{\it x3}{\it y1}^{2} }\\
  {\it PT^{x1}_{x1x2} }&=&\frac{1}{2}\frac{{\it y3} }{{\it x3}{\it y2} }\quad\quad\quad\quad  {\it PT^{x2}_{x1x2} }=\frac{3}{2}\frac{{\it y3} }{{\it x3}{\it y1} }\\
  {\it PT^{ x1}_{ x2 x2} }&=&-\frac{1}{2}\frac{{\it y1}{\it y3} }{{\it x3}{\it y2}^{2}}\quad\quad\quad  {\it PT^{x2}_{x2x2} }=-\frac{3}{2}\frac{{\it y3} }{{\it x3}{\it y2} }
 \end{eqnarray*}
\begin{Maple Normal}{
\begin{Maple Normal}{
{\textbf{Curvature tensors of Cartan  connection}}}\end{Maple Normal}
}\end{Maple Normal}

\vspace{.2cm}

\begin{maplegroup}
\begin{Maple Normal}{
$\bullet$ \textbf{ h-curvature   tensor ${R^h_{ijk}}$:}}\end{Maple Normal}
\end{maplegroup}

\begin{maplegroup}
\begin{mapleinput}
\mapleinline{active}{1d}{definetensor(RC[i,-h,-j,-k] = Hdiff(Gammastar[i,-h,-j],X[k])
-Hdiff(Gammastar[i,-h,-k], X[j])+Gammastar[m,-h,-j]
*Gammastar[i,-m, -k]-Gammastar[m,-h,-k]*Gammastar[i,-m,-j]
+C[i,-h,-m]*RG[m,-j,-k], antisymm[3,4]):}{\[\]}
\end{mapleinput}
\end{maplegroup}
\begin{maplegroup}
\begin{mapleinput}
\mapleinline{active}{1d}{show(RC[i,-h,-j,-k]); }{\[\]}
\end{mapleinput}
\mapleresult
\end{maplegroup}
{\vspace{-.5cm}}\begin{eqnarray*}
 {\it RC^{x1}_{x1x1x2} }&=&-\frac{3}{8}\frac{{\it y1}^{2}}{{\it x3}{\it y2}^{2} }\quad\quad\quad  {\it RC^{x1}_{x2x1x2} }=\frac{1}{4}\frac{{\it y1}^{3}}{{\it x3}{\it y2}^{3} }\\
  {\it RC^{x2}_{x1x1x2} }&=&-\frac{3}{4}\frac{{\it y1} }{{\it x3}{\it y2} }\quad\quad\quad  {\it RC^{x2}_{x2x1x2} }=\frac{3}{8}\frac{{\it y1}^{2}}{{\it x3}{\it y2}^{2} }\\
  {\it RC^{x2}_{x3x1x3} }&=&-\frac{1}{4{\it x3}^{2}}\quad\quad\quad\quad  {\it RC^{x3}_{x1x1x3} }=\frac{3}{4}\frac{{\it y1} }{{\it x3}{\it y2} }\\
   {\it RC^{x3}_{x2x1x3} }&=&-\frac{3}{8}\frac{{\it y1}^{2}}{{\it x3}{\it y2}^{2} }\quad\quad\quad  {\it RC^{x2}_{x3x2x3} }=-\frac{1}{4{\it x3}^{2}}\\
   {\it RC^{x3}_{x1x2x3} }&=&-\frac{3}{8}\frac{{\it y1}^{2}}{{\it x3}{\it y2}^{2} }\quad\quad\quad  {\it RC^{x3}_{x2x2x3} }=\frac{1}{4}\frac{{\it y1}^{3}}{{\it x3}{\it y2}^{3} }
\end{eqnarray*}
\begin{maplegroup}
\begin{Maple Normal}{
$\bullet$ \textbf{ hv-curvature   tensor $P^h_{ijk}$:}}\end{Maple Normal}
\end{maplegroup}

\begin{maplegroup}
\begin{mapleinput}
\mapleinline{active}{1d}{definetensor(FT[i,-j,-k,-h] = Hdiff(C[i,-j,-k], X[h])
+Gammastar[i,-h,-u]*C[u,-k,-j]-Gammastar[u,-k,-h]*C[i,-u,-j]
-Gammastar[u,-h,-j]*C[i,-u,-k]): }{\[\]}
\end{mapleinput}
\end{maplegroup}
\begin{maplegroup}
\begin{mapleinput}
\mapleinline{active}{1d}{definetensor(PC[i,-h,-j,-k] = tddiff(Gammastar[i,-h,-j],Y[k])
-FT[i,-h,-k,-j]+C[i,-h,-m]*PT[m,-j,-k]); }{\[\]}
\end{mapleinput}
\end{maplegroup}
{\vspace{-.5cm}}\begin{eqnarray*}
 {\it PC^{x1}_{x3x1x1} }&=&-\frac{1}{2{\it x3}{\it y1} }\quad\quad\quad  {\it PC^{ x1}_{ x3 x1 x2} }=\frac{1}{2{\it x3}{\it y2} }\\
 {\it PC^{ x1}_{ x3 x2 x1} }&=&\frac{1}{2{\it x3}{\it y2} }\quad\quad\quad\quad  {\it PC^{ x1}_{ x3 x2 x2} }=-\frac{1}{2}\frac{{\it y1} }{{\it x3}{\it y2}^{2}}\\
{\it PC^{ x2}_{ x3 x1 x1} }&=&-\frac{3}{2}\frac{{\it y2} }{{\it x3}{\it y1}^{2} }\quad\quad\,  {\it PC^{ x2}_{ x3 x1 x2} }=\frac{3}{2{\it x3}{\it y1} } \\
 {\it PC^{ x2}_{ x3 x2 x1} }&=&\frac{3}{2{\it x3}{\it y1} }\quad\quad\quad\quad  {\it PC^{ x2}_{ x3 x2 x2} }=-\frac{3}{2{\it x3}{\it y2} }\\
{\it PC^{ x3 }_{x1 x1 x1} }&=&-\frac{3}{4{\it y2} }\quad\quad\quad\quad  {\it PC^{ x3}_{ x1 x1 x2} }=\frac{3}{4}\frac{{\it y1} }{{\it y2}^{2}}\\
{\it PC^{ x3}_{ x1 x2 x1} }&=&\frac{3}{4}\frac{{\it y1} }{{\it y2}^{2}}\quad\quad\quad\quad {\it PC^{ x3}_{ x1 x2 x2} }=-\frac{3}{4}\frac{{\it y1}^{2}}{{\it y2}^{3}}\\
{\it PC^{ x3}_{ x2 x1 x1} }&=&\frac{3}{4}\frac{{\it y1} }{{\it y2}^{2}}\quad\quad\quad\quad  {\it PC^{ x3}_{ x2 x1 x2} }=-\frac{3}{4}\frac{{\it y1}^{2}}{{\it y2}^{3}}\\
{\it PC^{ x3}_{ x2 x2 x1} }&=&-\frac{3}{4}\frac{{\it y1}^{2}}{{\it y2}^{3}}\quad\quad\quad\,\,  {\it PC^{ x3}_{ x2 x2 x2} }=\frac{3}{4}\frac{{\it y1}^{3}}{{\it y2}^{4}}
\end{eqnarray*}
\begin{maplegroup}
\begin{Maple Normal}{
$\bullet$ \textbf{ v-curvature   tensor ${S^h_{ijk}}$:}}\end{Maple Normal}
\end{maplegroup}
\vspace{.2cm}
\begin{maplegroup}
\begin{mapleinput}
\mapleinline{active}{1d}{definetensor(S[i,-h,-j,-k] = C[m,-h,-k]*C[i,-m,-j]
-C[m,-h,-j]*C[i,-m,-k]):
}{}
\end{mapleinput}
\end{maplegroup}
\begin{maplegroup}
\begin{mapleinput}
\mapleinline{active}{1d}{show(S[i,-h,-j,-k]); }{\[\]}
\end{mapleinput}
\mapleresult
\begin{maplelatex}
\mapleinline{inert}{2d}{}{\[ S^{{h}}_{{\it i j k} }=0\]}
\end{maplelatex}
\end{maplegroup}

\vspace{.2cm}
\begin{Maple Normal}{
\begin{Maple Normal}{
{\textbf{Curvature tensors of Berwald connection}}}\end{Maple Normal}
}\end{Maple Normal}

\vspace{.2cm}

\begin{maplegroup}
\begin{Maple Normal}{
$\bullet$ \textbf{ h-curvature   tensor $\overcirc{R}^h_{ijk}$:}}\end{Maple Normal}
\end{maplegroup}
\vspace{.2cm}
\begin{maplegroup}
\begin{mapleinput}
\mapleinline{active}{1d}{definetensor(RB[i,-h,-j,-k]= Hdiff(G[i,-h,-j], X[k])
-Hdiff(G[i,-h,-k], X[j])+G[m,-h,-j]*G[i,-m,-k]
-G[m,-h,-k]*G[i,-m,-j], antisymm[3, 4]): }{\[\]}
\end{mapleinput}
\end{maplegroup}
\begin{maplegroup}
\begin{mapleinput}
\mapleinline{active}{1d}{show(RB[-i,h,-j,-k]);}{\[\]}
\end{mapleinput}
\end{maplegroup}
{\vspace{-.5cm}}\begin{eqnarray*}
  {\it RB^{x1}_{x1x1x2} }&=&-\frac{3}{8}\frac{{\it y1}^{2}}{{\it x3}{\it y2}^{2} }\quad\quad\quad  {\it RB^{x1}_{x2x1x2} }=\frac{1}{4}\frac{{\it y1}^{3}}{{\it x3}{\it y2}^{3} }\\
  \\
  {\it RB^{x2}_{x1x1x2} }&=&-\frac{3}{4}\frac{{\it y1} }{{\it x3}{\it y2} }\quad\quad\quad  {\it RB^{x2}_{x2x1x2} }=\frac{3}{8}\frac{{\it y1}^{2}}{{\it x3}{\it y2}^{2} }\\
  {\it RB^{x2}_{x3x1x3} }&=&-\frac{1}{4{\it x3}^{2}}\quad\quad\quad\quad  {\it RB^{x3}_{x1x1x3} }=\frac{3}{4}\frac{{\it y1} }{{\it x3}{\it y2} }\\
  {\it RB^{x3}_{x2x1x3} }&=&-\frac{3}{8}\frac{{\it y1}^{2}}{{\it x3}{\it y2}^{2} }\quad\quad\quad  {\it RB^{x2}_{x3x2x3} }=-\frac{1}{4{\it x3}^{2}}\\
  {\it RB^{x3}_{x1x2x3} }&=&-\frac{3}{8}\frac{{\it y1}^{2}}{{\it x3}{\it y2}^{2} }\quad\quad\quad  {\it RB^{x3}_{x2x2x3} }=\frac{1}{4}\frac{{\it y1}^{3}}{{\it x3}{\it y2}^{3} }
\end{eqnarray*}
\begin{maplegroup}
\begin{Maple Normal}{
$\bullet$ \textbf{ hv-curvature   tensor $\overcirc{P}^h_{ijk}$:}}\end{Maple Normal}
\end{maplegroup}

\vspace{.2cm}

\begin{maplegroup}
\begin{mapleinput}
\mapleinline{active}{1d}{definetensor(PB[i,-h,-j,-k]= tddiff(G[i,-h,-j],Y[k])): }{\[\]}
\end{mapleinput}
\end{maplegroup}
\begin{maplegroup}
\begin{mapleinput}
\mapleinline{active}{1d}{show(PB[h,-i,-j,-k]); }{\[\]}
\end{mapleinput}
\mapleresult
\end{maplegroup}
{\vspace{-.5cm}}\begin{eqnarray*}
  {\it PB^{ x3}_{ x1 x1 x1} }&=&-\frac{3}{2{\it y2} }\quad\quad\quad  {\it PB^{ x3}_{ x1 x1 x2} }=\frac{3}{2}\frac{{\it y1} }{{\it y2}^{2}}\\
 {\it PB^{ x3}_{ x1 x2 x2} }&=&-\frac{3}{2}\frac{{\it y1}^{2}}{{\it y2}^{3}}\quad\quad\quad  {\it PB^{ x3}_{ x2 x2 x2} }=\frac{3}{2}\frac{{\it y1}^{3}}{{\it y2}^{4}}
\end{eqnarray*}
\begin{Maple Normal}{
\begin{Maple Normal}{
{\textbf{Curvature tensors of Chern connection}}}\end{Maple Normal}
}\end{Maple Normal}

\vspace{.2cm}

\begin{maplegroup}
\begin{Maple Normal}{
$\bullet$ \textbf{ h-curvature   tensor $\overdiamond{R}^h_{ijk}$:}}\end{Maple Normal}
\end{maplegroup}

\vspace{.2cm}

\begin{maplegroup}
\begin{mapleinput}
\mapleinline{active}{1d}{definetensor(Rchern[i,-h,-j,-k] = Hdiff(Gammastar[i,-h,-j], X[k])
-Hdiff(Gammastar[i,-h,-k],X[j])+Gammastar[m,-h,-j]*Gammastar[i,-m,-k]
-Gammastar[m,-h,-k]*Gammastar[i,-m,-j], antisymm[3,4]): }{\[\]}
\end{mapleinput}
\end{maplegroup}
\begin{maplegroup}
\begin{mapleinput}
\mapleinline{active}{1d}{show(Rchern[i,-h,-j,-k]); }{\[\]}
\end{mapleinput}
\mapleresult
\end{maplegroup}
{\vspace{-.5cm}}\begin{eqnarray*}
  {\it Rchern^{x1}_{x1x1x2} }&=&-\frac{1}{8}\frac{{\it y1}^{2}}{{\it x3}{\it y2}^{2} }\quad\quad\quad  {\it Rchern^{x2}_{x2x1x2} }=-\frac{3}{8}\frac{{\it y1}^{2}}{{\it x3}{\it y2}^{2} }\\
  {\it Rchern^{ x1}_{ x1 x1 x3} }&=&-\frac{1}{4}\frac{{\it y3} }{{\it x3}^{2}{\it y1} }\quad\quad\quad  {\it Rchern^{ x1}_{ x2 x1 x3} }=\frac{1}{4}\frac{{\it y3} }{{\it x3}^{2}{\it y2} }\\
  {\it Rchern^{x2}_{x3x1x3} }&=&-\frac{1}{4{\it x3}^{2}}\quad\quad\quad\quad\,  {\it Rchern^{ x2}_{ x1 x1 x3} }=-\frac{3}{4}\frac{{\it y2}{\it y3} }{{\it x3}^{2}{\it y1}^{2}}\\
  {\it Rchern^{ x2 }_{x2 x1 x3} }&=&\frac{3}{4}\frac{{\it y3} }{{\it x3}^{2}{\it y1} }\quad\quad\quad  {\it Rchern^{x3}_{x1x1x3} }=\frac{3}{4}\frac{{\it y1} }{{\it x3}{\it y2} }\\
  {\it Rchern^{x3}_{x2x1x3} }&=&-\frac{3}{8}\frac{{\it y1}^{2}}{{\it x3}{\it y2}^{2} }\quad\quad\,  {\it Rchern^{ x1}_{ x1 x2 x3} }=\frac{1}{4}\frac{{\it y3} }{{\it x3}^{2}{\it y2} }\\
  {\it Rchern^{ x1}_{ x2 x2 x3} }&=&-\frac{1}{4}\frac{{\it y1}{\it y3} }{{\it x3}^{2}{\it y2}^{2}}\quad\,\,\, {\it Rchern^{ x2}_{ x1 x2 x3} }=\frac{3}{4}\frac{{\it y3} }{{\it x3}^{2}{\it y1} }\\
  {\it Rchern^{ x2}_{ x2 x2 x3} }&=&-\frac{3}{4}\frac{{\it y3} }{{\it x3}^{2}{\it y2} }\quad\quad\,  {\it Rchern^{x2}_{x3x2x3} }=-\frac{1}{4{\it x3}^{2}}\\
  {\it Rchern^{x3}_{x1x2x3} }&=&-\frac{3}{8}\frac{{\it y1}^{2}}{{\it x3}{\it y2}^{2} }\quad\quad\,  {\it Rchern^{x3}_{x2x2x3} }=\frac{1}{4}\frac{{\it y1}^{3}}{{\it x3}{\it y2}^{3} }
\end{eqnarray*}
\begin{maplegroup}
\begin{Maple Normal}{
$\bullet$ \textbf{ hv-curvature   tensor $\overdiamond{P}^h_{ijk}$:}}\end{Maple Normal}
\end{maplegroup}

\vspace{.2cm}

\begin{maplegroup}
\begin{mapleinput}
\mapleinline{active}{1d}{definetensor(Pchern[i,-h,-j,-k]=tddiff(Gammastar[i,-h,-j],Y[k])):}{\[\]}
\end{mapleinput}
\end{maplegroup}
\begin{maplegroup}
\begin{mapleinput}
\mapleinline{active}{1d}{show(Pchern[h,-i,-j,-k]);}{\[\]}
\end{mapleinput}
\mapleresult
\end{maplegroup}
{\vspace{-.5cm}}\begin{eqnarray*}
  {\it Pchern^{ x1}_{ x1 x1 x1} }&=&-\frac{1}{2}\frac{{\it y3} }{{\it x3}{\it y1}^{2} }\hspace{1.3cm}  {\it Pchern^{ x1}_{ x1 x1 x3} }=\frac{1}{2{\it x3}{\it y1} }\\
  {\it Pchern^{ x2 }_{x1 x1 x1} }&=&-\frac{3{\it y2}{\it y3} }{{\it x3}{\it y1}^{3} } \hspace{1.5cm}  {\it Pchern^{x2}_{x1x1x2} }=\frac{3}{2}\frac{{\it y3} }{{\it x3}{\it y1}^{2} }\\
  {\it Pchern^{ x2 }_{x1 x1 x3} }&=&\frac{3}{2}\frac{{\it y2} }{{\it x3}{\it y1}^{2} }\hspace{1.6cm} {\it Pchern^{ x3 }_{x1 x1 x1} }=-\frac{3}{2{\it y2} }\\
  {\it Pchern^{ x3 }_{x1 x1 x2} }&=&\frac{3}{2}\frac{{\it y1} }{{\it y2}^{2}} \hspace{2cm} {\it Pchern^{ x1}_{ x1 x2 x2} }=\frac{1}{2}\frac{{\it y3} }{{\it x3}{\it y2}^{2}}\\
  {\it Pchern^{ x1}_{ x1 x2 x3} }&=&-\frac{1}{2{\it x3}{\it y2} }\hspace{1.5cm}  {\it Pchern^{ x2 }_{x1 x2 x1} }=\frac{3}{2}\frac{{\it y3} }{{\it x3}{\it y1}^{2} }\\
  {\it Pchern^{ x2 }_{x1 x2 x3} }&=&-\frac{3}{2{\it x3}{\it y1} }\hspace{1.5cm}  {\it Pchern^{ x3}_{ x1 x2 x1} }=\frac{3}{2}\frac{{\it y1} }{{\it y2}^{2}}\\
  {\it Pchern^{ x3 }_{x1 x2 x2} }&=&-\frac{3}{2}\frac{{\it y1}^{2}}{{\it y2}^{3}}\hspace{1.7cm} {\it Pchern^{ x1 }_{x2 x2 x1} }=\frac{1}{2}\frac{{\it y3} }{{\it x3}{\it y2}^{2}}\\
  {\it Pchern^{ x1 }_{x2 x2 x2} }&=&-\frac{{\it y1}{\it y3} }{{\it x3}{\it y2}^{3}}\hspace{1.5cm}  {\it Pchern^{ x1 }_{x2 x2 x3} }=\frac{1}{2}\frac{{\it y1} }{{\it x3}{\it y2}^{2}}\\
  {\it Pchern^{ x2 }_{x2 x2 x2} }&=&-\frac{3}{2}\frac{{\it y3} }{{\it x3}{\it y2}^{2}}\hspace{1.2cm}  {\it Pchern^{ x2 }_{x2 x2 x3} }=\frac{3}{2{\it x3}{\it y2} }\\
  {\it Pchern^{ x3}_{ x2 x2 x1} }&=&-\frac{3}{2}\frac{{\it y1}^{2}}{{\it y2}^{3}}\hspace{1.6cm}  {\it Pchern^{ x3}_{ x2 x2 x2} }=\frac{3}{2}\frac{{\it y1}^{3}}{{\it y2}^{4}}
\end{eqnarray*}
\newpage
\begin{Maple Normal}{
\begin{Maple Normal}{
{\textbf{Curvature tensors of Hashiguchi connection}}}\end{Maple Normal}
}\end{Maple Normal}

\vspace{.2cm}

\begin{maplegroup}
\begin{Maple Normal}{
$\bullet$ \textbf{ h-curvature   tensor $\overast{R}^h_{ijk}$:}}\end{Maple Normal}
\end{maplegroup}

\begin{maplegroup}
\begin{mapleinput}
\mapleinline{active}{1d}{definetensor(RH[i,-h,-j,-k] = Hdiff(G[i,-h,-j], X[k])
-Hdiff(G[i,-h,-k],X[j])+G[m,-h,-j]*G[i,-m,-k]-G[m,-h,-k]*G[i,-m,-j]
+C[i,-h,-m]*RG[m,-j,-k], antisymm[3, 4]): }{\[\]}
\end{mapleinput}
\end{maplegroup}
\begin{maplegroup}
\begin{mapleinput}
\mapleinline{active}{1d}{show(RH[i,-h,-j,-k]);}{\[\]}
\end{mapleinput}
\mapleresult
\end{maplegroup}
\vspace{-.5cm}\begin{eqnarray*}
 {\it RH^{x1}_{x1x1x2} }&=&-\frac{5}{8}\frac{{\it y1}^{2}}{{\it x3}{\it y2}^{2} }\quad\quad\quad  {\it RH^{x1}_{x2x1x2} }=\frac{1}{2}\frac{{\it y1}^{3}}{{\it x3}{\it y2}^{3} }\\
  {\it RH^{x2}_{x1x1x2} }&=&-\frac{3}{2}\frac{{\it y1} }{{\it x3}{\it y2} }\quad\quad\quad\,\,  {\it RH^{x2}_{x2x1x2} }=\frac{9}{8}\frac{{\it y1}^{2}}{{\it x3}{\it y2}^{2} }
  \end{eqnarray*}
 \vspace{-.5cm}\begin{eqnarray*}
  {\it RH^{ x1}_{ x1 x1 x3} }&=&\frac{1}{4}\frac{{\it y3} }{{\it x3}{\it y1}^{2}}\quad\quad\quad  {\it RH^{ x1}_{ x2 x1 x3} }=-\frac{1}{4}\frac{{\it y3} }{{\it x3}^{2}{\it y2}}\\
  {\it RH^{x2}_{x3x1x3} }&=&-\frac{1}{4{\it x3}^{2}}\quad\quad\quad \, \, \,  {\it RH^{ x2}_{ x1 x1 x3} }=\frac{3}{4}\frac{{\it y2}{\it y3} }{{\it x3}^{2}{\it y1}^{2}}\\
   {\it RH^{ x2}_{ x2 x1 x3} }&=&-\frac{3}{4}\frac{{\it y3} }{{\it x3}{\it y1}^{2}}\quad\quad\, \,   {\it RH^{x3}_{x1x1x3} }=\frac{3}{4}\frac{{\it y1} }{{\it x3}{\it y2} }\\
   {\it RH^{x3}_{x2x1x3} }&=&-\frac{3}{8}\frac{{\it y1}^{2}}{{\it x3}{\it y2}^{2} }\quad\quad\, \,  {\it RH^{ x1 }_{x1 x2 x3} }=-\frac{1}{4}\frac{{\it y3} }{{\it x3}^{2}{\it y2}}\\
   {\it RH^{ x1}_{ x2 x2 x3} }&=&\frac{1}{4}\frac{{\it y1}{\it y3} }{{\it x3}^{2}{\it y2}^{2}}\quad\quad\quad  {\it RH^{ x2 }_{x1 x2 x3} }=-\frac{3}{4}\frac{{\it y3} }{{\it x3}{\it y1}^{2}}\\
    {\it RH^{ x2}_{ x2 x2 x3} }&=&\frac{3}{4}\frac{{\it y3} }{{\it x3}^{2}{\it y2}}\quad\quad\quad \, \, \, {\it RH^{x2}_{x3x2x3} }=-\frac{1}{4{\it x3}^{2}}\\
    {\it RH^{x3}_{x1x2x3} }&=&-\frac{3}{8}\frac{{\it y1}^{2}}{{\it x3}{\it y2}^{2} }\quad\quad\quad  {\it RH^{x3}_{x2x2x3} }=\frac{1}{4}\frac{{\it y1}^{3}}{{\it x3}{\it y2}^{3} }
\end{eqnarray*}
\begin{Maple Normal}{
\begin{Maple Normal}{
$\bullet$ \textbf{ hv-curvature   tensor $\overast{P}^h_{ijk}$:}}\end{Maple Normal}
}\end{Maple Normal}

\vspace{.2cm}

\begin{maplegroup}
\begin{mapleinput}
\mapleinline{active}{1d}{definetensor(PH1[i,-h,-k,-j]=Hdiff(C[i,-h,-k],X[j])+G[i,-m,-j]
*C[m,-h,-k]-G[m,-h,-j]*C[i,-m,-k]-G[m,-k,-j]*C[i,-h,-m]): }{\[\]}
\end{mapleinput}
\end{maplegroup}
\begin{maplegroup}
\begin{mapleinput}
\mapleinline{active}{1d}{definetensor(PH[i,-h,-j,-k] =
tddiff(G[i,-h,-j], Y[k])-PH1[i,-h,-k,-j], symm[2,4]); }{\[\]}
\end{mapleinput}
\end{maplegroup}
\begin{maplegroup}
\begin{mapleinput}
\mapleinline{active}{1d}{show(PH[h,-i,-j,-k]);}{\[\]}
\end{mapleinput}
\mapleresult
\end{maplegroup}
\vspace{-.5cm}\begin{eqnarray*}
  {\it PH^{ x1 }_{x1 x1 x1} }&=&\frac{1}{2}\frac{{\it y3} }{{\it x3}{\it y1}^{2} }\quad\quad\quad  {\it PH^{ x2}_{ x1 x1 x1} }=\frac{3{\it y2}{\it y3} }{{\it x3}{\it y1}^{3} }\\
  {\it PH^{ x2 }_{x1 x2 x1} }&=&-\frac{3}{2}\frac{{\it y3} }{{\it x3}{\it y1}^{2} }\quad\quad \, \,  {\it PH^{ x3 }_{x1 x1 x1} }=-\frac{3}{4{\it y2} }\\
  {\it PH^{ x3 }_{x1 x2 x1} }&=&\frac{3}{4}\frac{{\it y1} }{{\it y2}^{2}}\quad\quad\quad\quad \,   {\it PH^{ x1}_{ x1 x2 x2} }=-\frac{1}{2}\frac{{\it y3} }{{\it x3}{\it y2}^{2}}\\
  {\it PH^{x2}_{x1x1x2} }&=&-\frac{3}{2}\frac{{\it y3} }{{\it x3}{\it y1}^{2} }\quad\quad \, \,  {\it PH^{ x3}_{ x1 x1 x2} }=\frac{3}{4}\frac{{\it y1} }{{\it y2}^{2}}\\
  {\it PH^{ x3}_{ x1 x2 x2} }&=&-\frac{3}{4}\frac{{\it y1}^{2}}{{\it y2}^{3}}\quad\quad\quad\,\,  {\it PH^{ x1 }_{x1 x1 x3} }=-\frac{1}{2{\it x3}{\it y1} }\\
  {\it PH^{ x1}_{ x1 x2 x3} }&=&\frac{1}{2{\it x3}{\it y2} }\quad\quad\quad \, \,\, \, {\it PH^{ x2 }_{x1 x1 x3} }=-\frac{3}{2}\frac{{\it y2} }{{\it x3}{\it y1}^{2} }\\
  {\it PH^{ x2 }_{x1 x2 x3} }&=&\frac{3}{2{\it x3}{\it y1} }\quad\quad\quad \, \,\, \, {\it PH^{x1}_{x2x1x2} }=-\frac{1}{2}\frac{{\it y3} }{{\it x3}{\it y2}^{2}}\\
  {\it PH^{ x1 }_{x2 x2 x2} }&=&\frac{{\it y1}{\it y3} }{{\it x3}{\it y2}^{3}}\quad\quad\quad \, \,\, \, {\it PH^{ x2}_{ x2 x2 x2} }=\frac{3}{2}\frac{{\it y3} }{{\it x3}{\it y2}^{2}}\\
  {\it PH^{ x3}_{ x2 x1 x2} }&=&-\frac{3}{4}\frac{{\it y1}^{2}}{{\it y2}^{3}}\quad\quad\quad \, \, \, {\it PH^{ x3 }_{x2 x2 x2} }=\frac{3}{4}\frac{{\it y1}^{3}}{{\it y2}^{4}}\\
  {\it PH^{ x1 }_{x2 x1 x3} }&=&\frac{1}{2{\it x3}{\it y2} }\quad\quad\quad \, \,\, \, {\it PH^{ x1 }_{x2 x2 x3} }=-\frac{1}{2}\frac{{\it y1} }{{\it x3}{\it y2}^{2}}\\
  {\it PH^{ x2 }_{x2 x1 x3} }&=&\frac{3}{2{\it x3}{\it y1} }\quad\quad\quad \, \,\, \, {\it PH^{ x2 }_{x2 x2 x3} }=-\frac{3}{2{\it x3}{\it y2} }\\
\end{eqnarray*}
The v-curvature of Hashiguchi connection is the same as the v-curvature of Cartan connection.
\begin{rem}
\em{According to the above consideration, if we calculate the hv-curvature P of Cartan connection, in the example mentioned in Section 3, we find that the components $P^h_{ijk}$ vanish identically as expected.
}
\end{rem}

\Section{ Tensor simplification}
\par
It is well known that the simplification of  tensor expressions is not an easy task  \cite{Portugal3}. However, we have noted that if we have a complicated formula  of a geometric object, such as $P^h_{ijk}$, we can significantly simplify its expression as follows. We let the package compute the tensor $P_{hijk}:=g_{rh}P^r_{ijk}$  (instead of $P^h_{ijk}$) and ask it to show  the tensor $P^h_{ijk}$.

\smallskip
To illustrate this technique  let us consider the following example.

\smallskip
Let $M=\mathbb{R}^3$, $U=\{(x1,x2,x3;y1,y2,y3)\in \mathbb{R}^3 \times \mathbb{R}^3:  y1\neq 0,y2\neq 0,y3\neq 0\}$.  Let  $F$ be the Finsler structure  defined on $U$ by
$$F=(x1y2^3+y1^2y3)^{1/3}.$$

For example, let us  compute the component $S^1_{112}$ of the  v-curvature tensor $S^h_{ijk}$ of Cartan connection.
\smallskip
\begin{maplegroup}
\begin{mapleinput}
\mapleinline{active}{1d}{definetensor(SC[i,-h,-j,-k] = C[m,-h,-k]*C[i,-m,-j]
-C[m,-h,-j]*C[i,-m,-k]):
}{}
\end{mapleinput}
\mapleresult
\begin{maplelatex}
\mapleinline{inert}{2d}{}{\[ {\it SC}^{i }_{h j k }=C^{m }_{h k }C^{i }_{m j }-C^{m }_{hj }C^{i }_{mk }\]}
\end{maplelatex}
\end{maplegroup}
\begin{maplegroup}
\begin{mapleinput}
\mapleinline{active}{1d}{show(SC[i,-h,-j,-k]); }{\[\]}
\end{mapleinput}
\mapleresult
\begin{maplelatex}
\mapleinline{inert}{2d}{}{\[ {\it SC}^{x1 }_{x1 x1 x2 }=-\frac{1}{18}~\frac{{\it y3} {\it y1}\left(-x1{\it y2}^{3}+{\it y3}{\it y1}^{2}\right)x1{\it y2}^{2}\left({\it y3}{\it y1}^{2}-3x1{\it y2}^{3}\right)}{\left(x1{\it y2}^{3}+{\it y3}{\it y1}^{2}\right)^{4}}\]}
\mapleinline{inert}{2d}{}{\hspace{1.7cm}\[-\frac{2}{27}~\frac{{\it y3}^{2}{\it y1}^{3}x1{\it y2}^{2}\left({\it y3}{\it y1}^{2}+3~x1{\it y2}^{3}\right)}{\left(x1{\it y2}^{3}+{\it y3}{\it y1}^{2}\right)^{4}}\]}
\mapleinline{inert}{2d}{}{\hspace{1.7cm}\[+\frac{1}{36}\frac{{\it y1}x1{\it y2}^{2}\left(-x1{\it y2}^{3}+{\it y3}{\it y1}^{2}\right)~{\it y3}\left(-3~x1{\it y2}^{3}+5{\it y3}{\it y1}^{2}\right)}{\left(x1{\it y2}^{3}+{\it y3}{\it y1}^{2}\right)^{4}}\]}
\mapleinline{inert}{2d}{}{\hspace{1.7cm}\[+\frac{1}{54}~\frac{{\it y1}x1{\it y2}^{2}\left(4{\it y3}^{2}{\it y1}^{4}+21{\it y3}{\it y1}^{2}x1{\it y2}^{3}+9x1^{2}{\it y2}^{6}\right){\it y3} }{\left(x1{\it y2}^{3}+{\it y3}{\it y1}^{2}\right)^{4}}\]}
\end{maplelatex}
\end{maplegroup}

\smallskip
The above expression is complicated. But, in fact, if we lower the index $i$ in the above definition and use the command \emph{show(SC[i,-h,-j,-k])}, then we have the following simplification.
\smallskip
\begin{maplegroup}
\begin{mapleinput}
\mapleinline{active}{1d}{definetensor(SC[-i,-h,-j,-k] = C[m,-h,-k]*C[-i,-m,-j]
-C[m,-h,-j]*C[-i,-m,-k]):
}{}
\end{mapleinput}
\begin{maplelatex}
\mapleinline{inert}{2d}{}{\[ {\it SC}_{i h j k }=C^{m }_{h k }C_{i m j }-C^{m }_{hj }C_{i mk }\]}
\end{maplelatex}
\end{maplegroup}
\begin{maplegroup}
\begin{mapleinput}
\mapleinline{active}{1d}{show(SC[i,-h,-j,-k]); }{\[\]}
\end{mapleinput}
\mapleresult
\begin{maplelatex}
\mapleinline{inert}{2d}{}{\[ {\it SC}^{x1 }_{x1 x1 x2 }=\frac{1}{12}~\frac{{\it y3}{\it y1}x1{\it y2}^{2}}{\left(x1{\it y2}^{3}+{\it y3}{\it y1}^{2}\right)^{2}},\]}
\end{maplelatex}
\mapleresult
\end{maplegroup}
\smallskip
\noindent which is very simple compared with its expression before simplification.
\begin{rem}
\em{Be careful when you lower or raise an index, this index should be lowerable or raisable. For example, in the definition of $P^h_{ijk}$ we encounter the term $\paa_k \Gamma^i_{hj}$ (cf. \S 1). The index $i$ in this term can not be lowered since $g_{im}(\paa_k \Gamma^m_{hj})\neq \paa_k (g_{im}\Gamma^m_{hj})$. So we can not use the command  tddiff(Gammastar[-i,-h,-j], Y[k]). Such a problem can be treated as illustrated below:
 }
\end{rem}
\smallskip
\begin{maplegroup}
\begin{mapleinput}
\mapleinline{active}{1d}{definetensor(FT[i,-j,-k,-h]=Hdiff(C[i,-j,-k], X[h])
+Gammastar[i,-h,-u]*C[u,-k,-j]-Gammastar[u,-k,-h]*C[i,-u,-j]
-Gammastar[u,-h,-j]*C[i,-u,-k]); }{\[\]}
\end{mapleinput}
\end{maplegroup}
\begin{maplegroup}
\begin{mapleinput}
\mapleinline{active}{1d}{definetensor(PC[i,-h,-j,-k] = tddiff(Gammastar[i,-h,-j],Y[k])
-FT[i,-h,-k,-j]+C[i,-h,-m]*PT[m,-j,-k]); }{\[\]}
\end{mapleinput}
\mapleresult
\begin{maplelatex}
\mapleinline{inert}{2d}{}{\[ {\it PC}^{i }_{h j k }={\it tddiff} _{k }\left({\it Gammastar}^{i }_{h j }\right)-{\it FT}^{i }_{h k j }+C^{i }_{h m }~{\it PT}^{m }_{j k }\]}
\end{maplelatex}
\end{maplegroup}
\begin{maplegroup}
\begin{mapleinput}
\mapleinline{active}{1d}{show(PC[i,-h,-j,-k]);}{\[\]}
\end{mapleinput}
\mapleresult
\begin{maplelatex}
\mapleinline{inert}{2d}{}{\[ {\it PC}^{x1 }_{x1 x1 x1 }=\frac{1}{72}\frac{1}{\left({\it x1} {\it y2}^{3}+{\it y1}^{2}{\it y3} \right)^{2}{\it y1}}\Big(\frac{3{\it x1} {\it y2}^{6}\left(-{\it y1}^{2}{\it y3} +3{\it x1} {\it y2}^{3}\right){\it y1}^{2}{\it y3} }{\left({\it x1} {\it y2}^{3}+{\it y1}^{2}{\it y3} \right)^{2}}\]}
\mapleinline{inert}{2d}{}{\hspace{1.7cm}\[+\frac{10{\it y1}^{2}\left({\it y1}^{2}{\it y3} +3{\it x1} {\it y2}^{3}\right){\it y2}^{6}{\it y3} {\it x1} }{\left({\it x1} {\it y2}^{3}+{\it y1}^{2}{\it y3} \right)^{2}}+\frac{3}{2}\frac{{\it x1}^{2}{\it y2}^{9}\left(-{\it y1}^{2}{\it y3} +3{\it x1} {\it y2}^{3}\right)}{\left({\it x1} {\it y2}^{3}+{\it y1}^{2}{\it y3} \right)^{2}}\]}
\mapleinline{inert}{2d}{}{\hspace{1.7cm}\[+\frac{5\left({\it y1}^{2}{\it y3} +3{\it x1} {\it y2}^{3}\right){\it y2}^{9}{\it x1}^{2\\
\mbox{}}}{\left({\it x1} {\it y2}^{3}+{\it y1}^{2}{\it y3} \right)^{2}}+\frac{3}{2}\frac{{\it y2}^{3}{\it y1}^{4}\left(-{\it y1}^{2}{\it y3}+3{\it x1} {\it y2}^{3}\right){\it y3}^{2}}{\left({\it x1} {\it y2}^{3}+{\it y1}^{2}{\it y3} \right)^{2}}\]}
\mapleinline{inert}{2d}{}{\hspace{1.7cm}\[+\frac{5{\it y2}^{3}{\it y1}^{4}\left({\it y1}^{2}{\it y3} +3{\it x1} {\it y2}^{3}\right){\it y3}^{2}}{\left({\it x1} {\it y2}^{3}+{\it y1}^{2}{\it y3} \right)^{2}}-15{\it x1} {\it y2}^{6}+{\it y1}^{2}{\it y3} {\it y2}^{3}\Big)\]}
\end{maplelatex}
\end{maplegroup}
\smallskip
This component can be simplified using the above mentioned technique.
\smallskip
\begin{maplegroup}
\begin{mapleinput}
\mapleinline{active}{1d}{definetensor(FT[i,-j,-k,-h] = Hdiff(C[i,-j,-k], X[h])
+Gammastar[i,-h, u]*C[u,-k,-j]-Gammastar[u,-k,-h]*C[i,-u,-j]
-Gammastar[u,-h,-j]*C[i,-u,-k]);  }{\[\]}
\end{mapleinput}
\end{maplegroup}
\begin{maplegroup}
\begin{mapleinput}
\mapleinline{active}{1d}{definetensor(ST[i,-h,-j,-k] = tddiff(Gammastar[i,-h,-j], Y[k]));  }{\[\]}
\end{mapleinput}
\end{maplegroup}
\begin{maplegroup}
\begin{mapleinput}
\mapleinline{active}{1d}{definetensor(PC[-i,-h,-j,-k] = g[-m,-i]*ST[m,-h,-j,-k]
-FT[-i,-h,-k,-j]+C[-i,-h,-m]*PT[m,-j,-k]);  }{\[\]}
\end{mapleinput}
\mapleresult
\begin{maplelatex}
\mapleinline{inert}{2d}{}{\[ {\it PC}_{i h j k }=g _{m i }~{\it ST}^{m }_{h j k }-{\it FT}_{i h k j }+C_{i h m }~{\it PT}^{m }_{j k }\]}
\end{maplelatex}
\end{maplegroup}
\begin{maplegroup}
\begin{mapleinput}
\mapleinline{active}{1d}{show(PC[i,-h,-j,-k]);  }{\[\]}
\end{mapleinput}
\mapleresult
\begin{maplelatex}
\mapleinline{inert}{2d}{}{\[ {\it PC}^{x1 }_{x1 x1 x1 }=\frac{1}{16}\frac{{\it y2}^{3}}{y1\left(x1 {\it y2}^{3}+{\it y3}{\it y1}^{2}\right)},\]}
\end{maplelatex}
\end{maplegroup}
\smallskip
\noindent which is  simpler compared with its expression before simplification.


\Section{ Conclusion}
 In this paper, we have achieved four objectives concerning the FINSLER package   \cite{Rutz3}, \cite{hbfinsler}:

\smallskip
 $\bullet$ The wrong calculation of the components of the hv-curvature tensor $P^h_{ijk}$ of Cartan connection has been corrected

$\bullet$ Modifications have been made so that the h- and hv-curvatures of Cartan connection (and other geometric objects) could be computed in all dimensions (not only dimension $4$).

$\bullet$ The package has been extended to compute not only the  geometric objects associated with Cartan connection but also those associated with other fundamental connections of Finsler geometry. Other definitions can be added  similarly to the package.

$\bullet$ A technique for simplifying tensor expressions has been introduced.

\smallskip
Thanks to the FINSLER package, one is able to study various examples and counterexamples in Finsler and Riemannian geometries. For example, in \cite{ND-cartan} and \cite{ND-Zadeh}, we have studied interesting counterexamples in Finsler geometry.

\providecommand{\bysame}{\leavevmode\hbox
to3em{\hrulefill}\thinspace}
\providecommand{\MR}{\relax\ifhmode\unskip\space\fi MR }
\providecommand{\MRhref}[2]{%
  \href{http://www.ams.org/mathscinet-getitem?mr=#1}{#2}
} \providecommand{\href}[2]{#2}

\end{document}